\documentclass[12pt]{article}

\usepackage{amsmath,amssymb,amsthm}
\usepackage{mathtools}

\usepackage{indentfirst}

\usepackage{graphicx}
\usepackage{tikz}
\usetikzlibrary{arrows.meta,fit}
\usepackage[dvipsnames,svgnames,table]{xcolor}

\usepackage{enumitem}

\usepackage{hyperref}
\hypersetup{
    colorlinks,
    bookmarksnumbered=true
}

\newtheorem{theorem}{Theorem}
\newtheorem{lemma}[theorem]{Lemma}

\theoremstyle{definition}

\theoremstyle{remark}

\newcommand{\Sym}{\operatorname{Sym}}

\newcommand{\PP}{\mathbb{P}}

\newcommand{\abs}[1]{\left\lvert #1 \right\rvert}

\title{Random Tournaments Are Asymptotically Almost Surely Link-irregular}

\author{Lily Zhang\thanks{Department of Mathematics, University of California, Berkeley\@. Email address: {\tt pink@berkeley.edu}.}}
\date{October 17, 2026}

\begin{document}

\maketitle

\begin{abstract}
    \sloppy
    Bastien and Khormali (2025) conjectured that a nontrivial link-irregular tournament with $n$ vertices exists if and only if $n \ge 6$. Their work establishes the conjecture for $n \le 100$. We present a probabilistic argument for the existence of link-irregular tournaments with $n$ vertices when $n \ge 9$ to resolve the conjecture. In fact, we show that a random tournament is asymptotically almost surely link-irregular.
\end{abstract}

\vspace{3mm}



\section{Introduction}

The study of links of graphs was first introduced by Ali, Chartrand, and Zhang \cite{AliChartrandZhang}. This study has since been extended to links of directed graphs by Bastien and Khormali \cite{BastienKhormali}. However, for the practical purposes of this paper, we assert the definitions involving link for directed graphs only.

Given a vertex $v \in V(G)$ of a directed graph, we denote its neighborhood (set of in-neighbors and out-neighbors) by $N(v)$. The \textbf{link} of a vertex $v \in G$ is defined by
    \[ L(v) = G[N(v)]. \]
In other words, $L(v)$ is the subgraph formed by the neighborhood of $v \in G$. $G$ is \textbf{link-irregular} if for all distinct $u, v \in G$,
    \[ L(u) \not\cong L(v). \]

It's easy to see that link-irregular tournaments are directed graphs $T$ which satisfy for all distinct $u, v \in V(T)$:
    \[ T-u \not\cong T-v. \]
Bastien and Khormali have conjectured that nontrivial link-irregular tournaments with $n$ vertices exist if and only if $n \ge 6$. They prove that the conjecture holds up to $n = 100$ \cite{BastienKhormali}.

In what follows, we resolve the conjecture by proving that link-irregular tournaments with $n$ vertices exist when $n \ge 9$. Our work invokes a probabilistic argument, which does not provide explicit examples of such directed graphs. However, the robustness of probabilistic reasoning provides us with a proof that uniformly random tournaments with $n$ vertices are asymptotically almost surely link-irregular.

Both of the main results have been formalized in the Lean 4 proof assistant and are available on GitHub \cite{Lean4}.
\section{Preliminary Results}

Before beginning the proof, we introduce notation and facts that we will regularly use.

It goes without saying that we may assume without loss of generality that $V(T) = \{ 1, \dotsc, n \}$.

Let $T_n$ be a uniformly random tournament on $\{ 1, \dotsc, n \}$. For a set $S$ and nonnegative integer $k$, we use the notation
    \[ \binom Sk = \{ A \subseteq S: \abs{A} = k \} \]
to denote the set of subsets of $S$ with $k$ elements.

For a permutation $\sigma \in \Sym(S)$ of a finite set, the action $\underline{\cdot}: \langle \sigma \rangle \times \binom S2 \to \binom S2$ is defined by
    \[ \sigma^m \cdot \{ x,y \} = \{ \sigma^m(x), \sigma^m(y) \}. \]
We denote $o_2(\sigma) = \abs{\binom S2 / \langle \sigma \rangle}$ to be the amount of orbits of the action $\langle \sigma \rangle$ on $\binom S2$.

Let $g_n$ be the amount of simple undirected graphs with $n$ unlabeled vertices. The values of $g_n$ can be found in A000088, where we only use the values up to $n = 18$ in the proof.
\begin{lemma}[Values of $g_n$]\label{Lemma 1}
    \begin{align*}
        g_0  &= 1 \\
        g_1  &= 1 \\
        g_2  &= 2 \\
        g_3  &= 4 \\
        g_4  &= 11 \\
        g_5  &= 34 \\
        g_6  &= 156 \\
        g_7  &= 1,044 \\
        g_8  &= 12,346 \\
        g_9  &= 274,668 \\
        g_{10} &= 12,005,168 \\
        g_{11} &= 1,018,997,864 \\
        g_{12} &= 165,091,172,592 \\
        g_{13} &= 50,502,031,367,952 \\
        g_{14} &= 29,054,155,657,235,488 \\
        g_{15} &= 31,426,485,969,804,308,768 \\
        g_{16} &= 64,001,015,704,527,557,894,928 \\
        g_{17} &= 245,935,864,153,532,932,683,719,776 \\
        g_{18} &= 1,787,577,725,145,611,700,547,878,190,848.
    \end{align*}
\end{lemma}
\begin{proof}
    See A000088 \cite{OEIS}.
\end{proof}

\section{Proof of Existence}

Our goal is to prove that if $n \ge 9$, then there exists a link-irregular tournament with $n$ vertices, whence the conjecture follows from the facts we mentioned in the introduction. We do this via Theorem~\ref{Theorem 6}, which additionally shows that random tournaments are asymptotically almost surely link-irregular.

There is notation that we wish to introduce that's only applicable for Lemma~\ref{Lemma 2} and Theorem~\ref{Theorem 3}.

Given distinct $u,v \in \{ 1, \dotsc, n \}$ and a bijection
    \[ \varphi \colon \{ 1, \dotsc, n \} \setminus \{u\}
    \to \{ 1, \dotsc, n \} \setminus \{v\}, \]
define the event
    \[ B_\varphi=\{\varphi\text{ is an isomorphism from }T_n-u\text{ to }T_n-v\}. \]
For every $x \not= u$, draw the directed edge $x \longrightarrow \varphi(x)$. Using the convention that loops $x \longrightarrow x$ are directed cycles of length $1$, let $C_\varphi$ be the set of vertices on directed cycles, and define $\sigma_\varphi \colon C_\varphi \to C_\varphi$ via $\sigma_\varphi(x)=\varphi(x)$.

\begin{figure}[ht]
    \centering
    \begin{tikzpicture}[
        vertex/.style={circle,draw,minimum size=6mm,inner sep=0pt},
        cycle vertex/.style={vertex,fill=RoyalBlue!12},
        arrow/.style={-{Stealth[length=2mm]},thick}
    ]
        \node[font=\small\bfseries] at (2.1,2) {$x\longmapsto\varphi(x)$};

        \node[vertex] (v) at (0,1.15) {$v$};
        \node[vertex] (a) at (1.4,1.15) {$a$};
        \node[vertex] (b) at (2.8,1.15) {$b$};
        \node[vertex] (u) at (4.2,1.15) {$u$};
        \draw[arrow] (v) -- (a);
        \draw[arrow] (a) -- (b);
        \draw[arrow] (b) -- (u);

        \node[cycle vertex] (c1) at (0.45,-0.7) {$c_1$};
        \node[cycle vertex] (c2) at (1.75,-0.7) {$c_2$};
        \node[cycle vertex] (c3) at (1.1,-1.75) {$c_3$};
        \node[cycle vertex] (d) at (3.55,-1.25) {$d$};
        \draw[arrow] (c1) -- (c2);
        \draw[arrow] (c2) -- (c3);
        \draw[arrow] (c3) -- (c1);
        \draw[arrow] (d) edge[loop right] (d);
        \coordinate (loop clearance) at ([xshift=6mm]d.east);
        \node[
            draw=RoyalBlue,
            dashed,
            rounded corners,
            fit=(c1)(c2)(c3)(d)(loop clearance),
            inner sep=5pt,
            label=above:{\color{RoyalBlue}$C_\varphi$}
        ] {};
    \end{tikzpicture}
    \caption{A representative graph of $\varphi$. The highlighted vertices form $C_\varphi$.}
\end{figure}
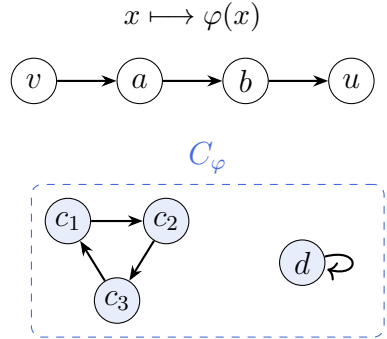

\begin{lemma}\label{Lemma 2}
    Let $n \ge 2$ be an integer, let $u,v \in \{1,\dotsc,n\}$ be distinct, and let\\
    $\varphi \colon \{ 1, \dotsc, n \} \setminus \{u\}
    \to \{ 1, \dotsc, n \} \setminus \{v\}$ be a bijection. Then
        \[ \PP(B_\varphi)\le 2^{-\binom{n-1}{2}+o_2(\sigma_\varphi)}. \]
\end{lemma}
\begin{proof}
Let $D_\varphi$ be the directed graph whose vertices are $V(D_\varphi) = \binom{\{ 1, \dotsc, n \}}2$, and the $\binom{n-1}2$ directed edges are
\[
    \{x,y\}\longrightarrow \{\varphi(x),\varphi(y)\}
\]
where $x,y \not= u$.

\begin{figure}[ht]
    \centering
    \begin{tikzpicture}[
        pair/.style={draw,rounded corners,fill=black!4,inner sep=3pt,font=\scriptsize},
        cycle pair/.style={pair,fill=RoyalBlue!12},
        arrow/.style={-{Stealth[length=2mm]},thick}
    ]
        \node[pair] (p1) at (0,1.15) {$\{v,c_1\}$};
        \node[pair] (p2) at (1.8,1.15) {$\{a,c_2\}$};
        \node[pair] (p3) at (3.6,1.15) {$\{b,c_3\}$};
        \node[pair] (p4) at (5.4,1.15) {$\{u,c_1\}$};
        \draw[arrow] (p1) -- (p2);
        \draw[arrow] (p2) -- (p3);
        \draw[arrow] (p3) -- (p4);

        \node[cycle pair] (q1) at (0.8,-0.15) {$\{c_1,c_2\}$};
        \node[cycle pair] (q2) at (2.7,-1.15) {$\{c_2,c_3\}$};
        \node[cycle pair] (q3) at (4.6,-0.15) {$\{c_3,c_1\}$};
        \draw[arrow] (q1) -- (q2);
        \draw[arrow] (q2) -- (q3);
        \draw[arrow] (q3) to[bend left=22] (q1);
    \end{tikzpicture}
    \caption{A path component and a cycle component of $D_\varphi$.}
\end{figure}
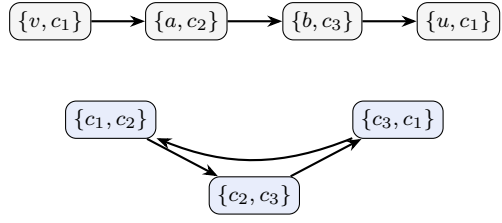

A vertex $\{x,y\} \in V(D_\varphi)$ lies on a directed cycle if and only if $x,y \in C_\varphi$, indicating that the directed cycles of $D_\varphi$ are precisely the orbits of $f: \binom{C_\varphi}2 \to \binom{C_\varphi}2$ by
\[
    f(\{x,y\}) = \{\sigma_\varphi(x),\sigma_\varphi(y)\}.
\]
So there are $o_2(\sigma_\varphi)$ connected components that are directed cycles.

Since $\varphi$ is a bijection, each vertex $\{x,y\} \in V(D_\varphi)$ satisfies\\
$\deg^-(\{x,y\}), \deg^+(\{x,y\}) \le 1$. So every connected component of $D_\varphi$ is either a directed path or a directed cycle. Therefore the amount of connected components of $D_\varphi$ is
    \[ r = \abs{V(D_\varphi)} - \abs{E(D_\varphi)} + o_2(\sigma_\varphi) = \binom n2 - \binom{n-1}2 + o_2(\sigma_\varphi). \]

Assuming $B_\varphi$ occurs, the orientation of one unordered pair in $T_n$ determines the orientations of all unordered pairs within this connected component. We conclude at most $2^r$ tournaments for which $B_\varphi$ occurs, thus
    \[ \PP(B_\varphi) \le \frac{2^r}{2^{\binom n2}} = \frac{2^{\binom n2 - \binom{n-1}2 + o_2(\sigma_\varphi)}}{2^{\binom n2}} = 2^{-\binom{n-1}{2}+o_2(\sigma_\varphi)}, \]
exactly as desired.
\end{proof}

\begin{theorem}\label{Theorem 3}
    Let $n \ge 2$ be an integer. Let $u, v \in \{ 1, \dotsc, n \}$ be distinct. Then
        \[ \PP(T_n-u \cong T_n-v) \le (n-2)! \cdot 2^{-\binom{n-1}2} \sum_{k=0}^{n-2} g_k. \]
\end{theorem}
\begin{proof}
    Taking the union over all bijections\\
    $\varphi: \{ 1, \dotsc, n \} \setminus \{u\} \to \{ 1, \dotsc, n \} \setminus \{v\}$ gives the equality of events
        \[ \{T_n-u \cong T_n-v\} = \bigcup_\varphi B_\varphi. \]
    Applying the union bound and Lemma~\ref{Lemma 2}, we obtain
        \[ \PP(T_n-u \cong T_n-v) \le \sum_\varphi \PP(B_\varphi) \le 2^{-\binom{n-1}2} \sum_\varphi 2^{o_2(\sigma_\varphi)}. \]
    Our goal is to prove the sufficient inequality
        \[ \sum_\varphi 2^{o_2(\sigma_\varphi)}
        \le (n-2)! \cdot \sum_{k=0}^{n-2}g_k. \]
    Fix $k \in \{ 0, \dotsc, n-2 \}$. We wish to evaluate the left sum restricted to $\abs{C_\varphi} = k$
        \[ \sum_{\substack{\varphi \\ \abs{C_\varphi}=k}} 2^{o_2(\sigma_\varphi)}. \]
    \begin{itemize}    
    \item There are $\binom{n-2}k$ ways to choose the $k$ vertices out of the $n-2$ vertices that aren't $u, v$ to determine $C_\varphi$.

    \item For a fixed $C\in\binom{\{1,\dotsc,n\}\setminus\{u,v\}}k$ and permutation $\sigma:C\to C$, there are $(n-k-2)!$ ways to order the remaining $n-k-2$ vertices to form the directed path from $v$ to $u$. So there are $(n-k-2)!$ bijections $\varphi$ with $C_\varphi=C$, $\sigma_\varphi=\sigma$.
    
    \item For a fixed $C_\varphi$, summing $2^{o_2(\sigma_\varphi)}$ over all permutations of $C_\varphi$ yields
    \[ \sum_{\sigma \in \Sym(C_\varphi)} 2^{o_2(\sigma)} = \sum_{\sigma\in S_k} 2^{o_2(\sigma)}. \]
\end{itemize}
Taking the sum from $k=0$ to $n-2$, we obtain
\begin{align*}
    \sum_\varphi 2^{o_2(\sigma_\varphi)}
    &= \sum_{k=0}^{n-2} \sum_{\substack{\varphi \\ \abs{C_\varphi}=k}} 2^{o_2(\sigma_\varphi)} \\
    &= \sum_{k=0}^{n-2} \binom{n-2}k \cdot (n-k-2)! \cdot \sum_{\sigma \in S_k} 2^{o_2(\sigma)} \\
    &= \sum_{k=0}^{n-2} \binom{n-2}k \cdot (n-k-2)! \cdot k! \cdot g_k \\
    &= (n-2)! \cdot \sum_{k=0}^{n-2} g_k,
\end{align*}
where the third equality follows from Burnside's Lemma. In fact, the sufficient inequality we wanted to prove is an equality!
\end{proof}

\begin{lemma}\label{Lemma 4}
    Let $n \ge 8$ be an integer. Then
        \[ g_n < 2 \cdot \frac{2^{\binom{n}2}}{n!}. \]
\end{lemma}
\begin{proof}
    Buhrman, Li, Tromp, and Vitányi \cite[Corollary~2]{BuhrmanLiTrompVitanyi} already supply the asymptotically stronger
        \[ g_n \le \frac{2^{\binom{n}2}}{n!} \left( 1 + \frac{4n^4}{2^n} \right). \]
    We merely show that $2^n > 4n^4$ for $n \ge 19$ via induction, then we perform a finite check on $n \in [8, 18]$.

    Base case ($n=19$): $2^{19} = 524288 > 521284 = 4 \cdot 19^4$.

    Induction step: Assuming it for $n-1$, we just show after dividing by $2^{n-1} > 4(n-1)^4$ that
        \[ 2 > \frac{n^4}{(n-1)^4}. \]
    This is true since the fourth root of the right side is
        \[ \frac n{n-1} \le \frac{20}{19} < 1.18 < \sqrt[4]2. \]
    
    Now we establish the inequality for $n \in [8, 18]$. This is shown through exhaustive casework using the values of $g_n$ explicitly listed in Lemma~\ref{Lemma 1}. Let $h_n = 2 \cdot \frac{2^{\binom{n}2}}{n!}$, so we need to show $g_n < h_n$. To not take up so much space, I present the values of $g_n$ and rounded-down integer estimates of $h_n$ in the following table. These numbers were obtained by a simple C++ program.
    \begin{center}
        \resizebox{1\textwidth}{!}{%
        \begin{tabular}{c|c|c}
            $n$ & $g_n$ & $h_n$ (rounded down) \\ \hline
            $8$  & $12{,}346$ & $13{,}315$ \\
            $9$  & $274{,}668$ & $378{,}744$ \\
            $10$ & $12{,}005{,}168$ & $19{,}391{,}739$ \\
            $11$ & $1{,}018{,}997{,}864$ & $1{,}805{,}194{,}655$ \\
            $12$ & $165{,}091{,}172{,}592$ & $308{,}086{,}554{,}595$ \\
            $13$ & $50{,}502{,}031{,}367{,}952$ & $97{,}070{,}963{,}663{,}284$ \\
            $14$ & $29{,}054{,}155{,}657{,}235{,}488$ & $56{,}800{,}381{,}023{,}544{,}553$ \\
            $15$ & $31{,}426{,}485{,}969{,}804{,}308{,}768$ & $62{,}041{,}162{,}845{,}983{,}597{,}115$ \\
            $16$ & $64{,}001{,}015{,}704{,}527{,}557{,}894{,}928$ & $127{,}060{,}301{,}508{,}574{,}406{,}891{,}621$ \\
            $17$ & $245{,}935{,}864{,}153{,}532{,}932{,}683{,}719{,}776$ & $489{,}824{,}936{,}450{,}937{,}195{,}885{,}253{,}015$ \\
            $18$ & $1{,}787{,}577{,}725{,}145{,}611{,}700{,}547{,}878{,}190{,}848$ & $3{,}566{,}796{,}337{,}249{,}846{,}674{,}392{,}882{,}402{,}392$
        \end{tabular}%
        }
    \end{center}
    
\end{proof}

\begin{lemma}\label{Lemma 5}
    Let $n \ge 7$ be an integer. Then
        \[ \sum_{k=0}^n g_k < 4 \cdot \frac{2^{\binom{n}2}}{n!} \]
\end{lemma}
\begin{proof}
    We show this via induction on $n$.

    Base case ($n = 7$): Just compare
    \begin{align*}
        g_0 + \dotsb + g_7 = 1+1+2&+4+11+34+156+1044 = 1253. \\
        4 \cdot \frac{2^{\binom{7}2}}{7!} &= \frac{2^{23}}{5040} = 1664 + \frac{128}{315}.
    \end{align*}

    Induction step: Assuming it true for $n-1$, we derive using Lemma~\ref{Lemma 4} that
        \[ \sum_{k=0}^n g_k < 4 \cdot \frac{2^{\binom{n-1}2}}{(n-1)!} + 2 \cdot \frac{2^{\binom{n}2}}{n!}. \]
    It's enough to show
        \[ 2 \cdot \frac{2^{\binom{n-1}2}}{(n-1)!} \le \frac{2^{\binom{n}2}}{n!}, \]
    where division reduces this to
        \[ \frac{2^{\frac{n(n-1) - (n-1)(n-2)}2}}{n} = \frac{2^{n-1}}{n} \ge 2. \]
    Obviously this is true for $n \ge 8$.
\end{proof}

\begin{theorem}\label{Theorem 6}
    If $n \ge 9$ is an integer, then there exists a link-irregular tournament with $n$ vertices. In fact, $T_n$ is asymptotically almost surely link-irregular.
\end{theorem}
\begin{proof}
    Let $u, v \in \{ 1, \dotsc, n \}$ be distinct.
    
    If $n-2 \ge 8$, combining Theorem~\ref{Theorem 3} and Lemma~\ref{Lemma 5} results in the chain of inequalities
    \begin{align*}
        \PP(T_n-u \cong T_n-v)
        &\le (n-2)! \cdot 2^{-\binom{n-1}2} \sum_{k=0}^{n-2} g_k \\
        &< (n-2)! \cdot 2^{-\binom{n-1}2} \cdot 4 \cdot \frac{2^{\binom{n-2}2}}{(n-2)!} \\
        &= 4 \cdot 2^{\frac{(n-2)(n-3) - (n-1)(n-2)}{2}} \\
        &= 2^{-n + 4}.
    \end{align*}
    Taking the union bound over all $\binom{n}2$ unordered pairs $u \not= v$ of $\{ 1, \dotsc, n \}$, we derive
        \[ \PP(T_n \text{ is not link-irregular}) < \underbrace{\binom{n}2 \cdot 2^{-n + 4}}_{a_n} \; \overset{n \to \infty}{\longrightarrow} \; 0. \]
    This gives us the second part of the theorem.

    It's not hard to show that $a_n$ is decreasing in $n$, and also
        \[ a_{10} = \binom{10}2 \cdot 2^{-10 + 4} = \frac{45}{64} < 1. \]
    Hence the first part of the theorem for $n \ge 10$.

    We handle the $n=9$ case separately, which follows via directly applying a similar union bound and Theorem~\ref{Theorem 3}:
    \begin{align*}
        \PP(T_9 \text{ is not link-irregular})
        &\le \binom92 \cdot 7! \cdot 2^{-\binom82}(g_0 + \dotsb + g_7) \\
        &= \frac{36 \cdot 7! \cdot 1253}{2^{28}} \\
        &= \frac{3552255}{4194304} \\
        &< 1.
    \end{align*}
\end{proof}

\bibliographystyle{plain}
\bibliography{references}

\end{document}